\documentclass[12pt,reqno]{amsart}
\usepackage{amsmath, amssymb, amsfonts, amsthm}
\usepackage{mathtools,xcolor}
\allowdisplaybreaks
\usepackage{hyperref}
\usepackage[all,cmtip]{xy}
\usepackage{tkz-euclide}
\usepackage{tikz-cd}
\usepackage[most]{tcolorbox}   
\usepackage{xcolor}            
\newtheorem{thm}{Theorem}[section]
\newtheorem*{thm*}{Theorem}
\newtheorem{lem}[thm]{Lemma}
\newtheorem{prop}[thm]{Proposition}
\theoremstyle{definition}
\newtheorem{dfn}[thm]{Definition}

\theoremstyle{plain}
\newtheorem{cor}[thm]{Corollary}
\numberwithin{equation}{section}
\usepackage[a4paper, top = 1.2in, bottom = 1.2in, left = 1.2in, right = 1.2in]{geometry}
\numberwithin{equation}{section}

\newcommand{\A}{\mathbb{A}}
\newcommand{\C}{\mathbb{C}}
\newcommand{\N}{\mathbb{N}}
\newcommand{\Q}{\mathbb{Q}}

\newcommand{\Z}{\mathbb{Z}}

\newcommand{\F}{\mathbb{F}}

\newcommand{\mcO}{\mathcal{O}}
\newcommand{\mcS}{\mathcal{S}}

\newcommand{\mfb}{\mathfrak{b}}

\newcommand{\mfp}{\mathfrak{p}}

\newcommand{\mfn}{\mathfrak{n}}

\newcommand{\SL}{\mathrm{SL}}
\newcommand{\GL}{\mathrm{GL}}

\newcommand{\Aut}{\mathrm{Aut}}

\newcommand{\Frob}{\mathrm{Frob}}

\newcommand{\Spec}{\mrm{Spec}}

\newcommand{\trace}{\mrm{trace}}

\newcommand{\Gal}{\mrm{Gal}}

\newcommand{\ra}{\rightarrow}

\newcommand{\mrm}[1]{\mathrm{#1}}

\def\1{1\!\!1}

\title[Coefficient fields of primitive HMFs and Nekov{\'a}{\v r}'s Theorem]{On generation of the coefficient field of a primitive HMF by a single Fourier coefficient: An application of Nekov{\'a}{\v r}'s Theorem}
\author[N. Kumar]{Narasimha Kumar}
\email{narasimha@math.iith.ac.in}
\address{
	Department of Mathematics \\
	Indian Institute of Technology Hyderabad\\
	Kandi, Sangareddy - 502284\\
	INDIA.
}
\author[D. Shit]{Dwipanjana Shit}
\email{ma22resch01001@iith.ac.in}
\address{
	Department of Mathematics \\
	Indian Institute of Technology Hyderabad\\
	Kandi, Sangareddy - 502284\\
	INDIA.
}

\begin{document}
 
    \begin{abstract}
         Let $f$ be a non-CM primitive Hilbert modular form (HMF) over a totally real number field, of parallel even weight $k\geq 2$, with coefficient field $E_f$, and let $F_f \subseteq E_f$ denote the field fixed by the inner twists of $f$. 
         In this article, we establish an analogue of Ribet’s large image theorem for the residual Galois representations associated with $ f $ (cf. Theorem~\ref{analogue_of_Ribet}), by applying  Nekov{\'a}{\v r}'s large image theorem (cf.~\cite[Theorem B.5.2]{Nek12}) to the automorphic representation $\pi(f)$ associated with $f$. As a result, we prove that if $f$ has no nontrivial inner twists, then $E_f$ is generated by a single prime Fourier coefficient, for a set of prime ideals of Dirichlet density one. We derive this result as a consequence of a more general theorem established for $F_f$. 
    \end{abstract}
	\maketitle

\section{Introduction and Statement of the Main Theorem}

    The arithmetic of Fourier coefficients of a primitive form is a classical theme in number theory, as these coefficients encode rich arithmetic information and satisfy strong algebraic relations through the action of Hecke operators. They also provide a fundamental link between classical arithmetic problems, $L$-functions, Galois representations, and arithmetic geometry.

    Let $f=\sum_{n=1}^{\infty}a_n(f)q^n\in S_k(\Gamma_1(N),\chi)^{\mathrm{new}}$ be a primitive form, with coefficient field $E_f$, and let $F_f$ denote the subfield of $E_f$ fixed by the inner twists of $f$. Maeda's conjecture predicts that for every positive integer $n$, the characteristic polynomial of the Hecke operator $T_n$ acting on $S_k(\SL_2(\Z))$, is irreducible over $\Q$
    and has Galois group $S_d$, where $d=\dim S_{k}(\SL_2(\Z))$ (cf.~\cite{BM03}). We refer to these characteristic polynomials as Hecke polynomials. In particular,  the space $S_k(\SL_2(\Z))$ consists of a single Galois conjugacy class.

    For higher levels $N > 1$, the generalized Maeda's conjecture predicts that, after decomposing the newspace into its natural arithmetic components, each component consists of a single Galois conjugacy class of newforms, denoted by $[f]$. Under this framework, for any prime $p$, the characteristic polynomial of $T_p$ on that component is irreducible over $\Q$ if and only if the degree $[E_f: \Q]$ equals the dimension of the corresponding component, and  for any prime $p\nmid N$, one has $E_f= \Q(a_p(f))$. This leads to the question of whether and for how many primes $p$ the coefficient field $E_f$ is generated by $a_p(f)$.

    \subsection*{Literature Survey}
    Let $\mathbb{P}$ denote the set of all prime numbers. Initial progress on the irreducibility of Hecke polynomials, acting on $S_k(\SL_2(\Z)),$  was made by James and Ono \cite{JO98}; they proved that a positive proportion of Hecke polynomials of $T_p\ (p \in \mathbb{P})$ are irreducible over $\mathbb{Q}$, under some assumption on Hecke polynomials. In \cite{BM03}, Baba and Murty significantly deepened this connection by showing that if the Hecke polynomial of $T_q$ is irreducible for some $q \in \mathbb{P}$, then the Hecke polynomials of $T_p$ remain irreducible for a density one set of primes $p$, establishing that $E_f = \mathbb{Q}(a_p(f))$ holds for a set of primes $p$ of density one.

    At higher levels $N>1$, the presence of inner twists obstructs the generation of $E_f$ by a single prime Fourier coefficient. In \cite{KSW08}, Koo, Stein, and Wiese resolved this question in the setting of classical primitive forms across arbitrary levels and characters. Using Ribet's large image theorem (cf.~\cite[Theorem 3.1]{Rib85}, \cite[Theorem 3.1]{KSW08}) for the residual Galois representations associated with $f$, together with the Chebotarev density theorem, they proved that $\Q(a_p(f)^2/\chi(p))= F_f$, for a density one set of primes $p$. Moreover, if $f$ has no nontrivial inner twists, they proved that $E_f=\mathbb{Q}(a_p(f))$ holds for a density one set of primes $p$, without requiring any prior assumptions on the irreducibility of Hecke polynomials.

    Naturally, one is led to the analogous question for primitive Hilbert modular forms $f$ over a totally real number field $F$ of parallel weight $k$. In \cite{KS23}, the first author and Sahoo studied this problem and established density one results for the generation of both $E_f$ and $F_f$, under the assumption that an analogue of Ribet’s large image theorem holds for the residual Galois representations attached to $f$.
    However, to establish a version of this large image theorem for $f$ of weight $k=2$ (cf.~\cite[Proposition 3.3]{KS23}), they had to work with the restrictive condition $[F:\Q]=[E_f:\Q]$ being an odd prime. This analysis relied on the work of Dimitrov \cite[Propositions 3.8 and 3.9]{Dim05} and Dieulefait–Dimitrov \cite{DD06}. Finally, they also provided explicit examples of primitive Hilbert modular forms of weight $2$ satisfying this restrictive hypothesis on degrees.

\subsection{Our contribution:}
    Let $f$ be a primitive Hilbert modular forms   over a totally real number field $F$ of  parallel even weight $k \geq 2$. In this work, we establish that an analogue of Ribet’s large image theorem holds for the residual Galois representations attached to $f$ (cf. Theorem~\ref{analogue_of_Ribet} in the text).  This objective is achieved from the properties of the automorphic representation $\pi(f)$ associated with $f$. More precisely, to obtain an analogue of Ribet's large image theorem for the Galois representations attached to $f$, we use the Nekov{\'a}{\v r}'s large image theorem~\cite[Theorem~B.5.2]{Nek12} for $\pi(f)$. This passage requires a careful comparison of the coefficient fields, the inner twists, and characters of $\pi(f)$ with those of $f$. Consequently, we establish that if $f$ has no nontrivial inner twists, then $E_f$ is generated by a single prime Fourier coefficient, for a set of prime ideals of Dirichlet density one; in fact, we derive this result as a consequence of a more general theorem established for $F_f$. An analogous question can be studied for intermediate subfields with the corresponding densities being completely determined by Corollary~\ref{density_generation_intermediate}.  We now state the main theorem of this article.

  \begin{thm}[Main Theorem]
    \label{main_thm}
	    Let $f$ be a non-CM primitive Hilbert modular form over a totally real number field $F$, of parallel even weight $k\geq 2$, level $\mfn$, and character $\Psi$, where $\Psi$ is as in \S \ref{HMF-Subsection}. Let $E_f:=\Q(C(\mfb,f) \mid (\mfb,\mfn)=1)$ be the field of coefficients of $f$ and $F_f:=\Q\big(\frac{C(\mfb,f)^2}{\Psi^{\ast}(\mfb)}\mid (\mfb,\mfn)=1\big)$ where $\Psi^{\ast}$ is the ideal character corresponding to $\Psi$. Then the set
        \begin{equation*}
            \left\{ \mfp \in \Spec(\mcO_F) \;\middle|\; \Q\left(\frac{C(\mfp,f)^2}{\Psi^{\ast}(\mfp)}\right) = F_f  \right\}
        \end{equation*}
        has Dirichlet density $1$. Furthermore, if $f$ does not have nontrivial inner twists, then the set
            $\left\{ \mfp \in \Spec(\mcO_F) \;\middle|\; \Q(C(\mfp, f)) = E_f \right\}$
        has Dirichlet density $1$.
	\end{thm}

\section{Preliminaries}
    In this section, we recall the necessary background on Hilbert modular forms, automorphic representations, Galois representations, and inner twists. Throughout this paper, we use the following notations:
\begin{itemize}
    \item Let $\mathbb{P}$ denote the set of all rational prime numbers.
    \item Let $F$ be a totally real number field of degree $n = [F : \mathbb{Q}]$, with ring of integers $\mcO_F$ and $P = \Spec(\mcO_F) \setminus \{(0) \}$. Let $\mfn \subseteq \mcO_F$ be a nonzero ideal.
    \item For any $\mfp\in P$, let $N(\mfp):=\#(\mcO_F/\mfp)$ denote its absolute norm. For a set $S\subseteq P$, the Dirichlet density of $S$ is defined by $\delta_D(S):=\lim_{s\to1^+}\frac{\sum_{\mfp\in S}N(\mfp)^{-s}}{\sum_{\mfp\in P}N(\mfp)^{-s}}$ if the limit exists.
    \item Let $k = (k_1, \ldots, k_n) \in \Z^n$ with $k_i \ge 2$ such that $k_1\equiv k_2\equiv\cdots\equiv k_n\pmod{2}$.
    Set $k_{0}:=\max\{k_1, \ldots, k_n\}$. We say that the weight $k$ is parallel if $k_1=\cdots=k_n$; in this case, we identify the weight vector $k$ with  $k_1$.
    \item For any extension $F \subseteq K  \subseteq \overline F$, set $G_K:=\Gal(\overline F/K).$ In particular, $G_F=\Gal(\overline F/F)$ is the absolute Galois group of $F$.
    \item For every unramified $\mfp\in P$, we denote the geometric Frobenius at $\mfp$ by $\Frob_{\mfp}$.
\end{itemize}

\subsection{Hilbert Modular Forms}
\label{HMF-Subsection}
    Let $\Psi:\A_{F}^{\times}/F^{\times}\rightarrow \C^{\times}$ be a Hecke character whose conductor divides $\mfn$ and whose infinity type is $2-k_0$. We denote by $S_{k}(\mfn, \Psi)$ the space of all Hilbert modular newforms over $F$ of weight $k$, level $\mfn$ and character $\Psi$. An eigenform $f \in S_{k}(\mfn, \Psi)$ is called a primitive form if it is normalised.

    Let $f \in S_k(\mathfrak{n}, \Psi)$ be a primitive form, whose Fourier coefficients $C(\mathfrak{b}, f)$ are parameterized by integral ideals $\mfb \subseteq \mcO_F$.
    For $\nu\in P$ with $\nu \nmid \mfn$, these coefficients coincide with the Hecke eigenvalues, i.e., $ T(\nu)f = \lambda_f(\nu)f$, where $\lambda_f(\nu) = C(\nu, f)$.

\subsubsection{Galois representations associated with primitive forms}
    Let $f \in S_{k}(\mfn, \Psi)$ be a primitive form with coefficient field $E_f$. By the works of Ohta, Carayol, Blasius-Rogawski, and Taylor (cf.~\cite{Car86}, \cite{Tay89}), for every $\lambda \in \Spec(\mcO_{E_f}) $, there exists a continuous $\lambda$-adic Galois representation
    \begin{equation*}
        \rho_{f, \lambda}:G_F \longrightarrow \Aut_{E_{f, \lambda}}\bigl(V_{\lambda}(f)\bigr)\cong\GL_2(E_{f, \lambda}),
    \end{equation*}
    where $E_{f, \lambda}$ denotes the completion of $E_f$ at $\lambda$, and $V_{\lambda}(f)$ is a two-dimensional $E_{f, \lambda}$-vector space (cf.~\cite[\S 0.4]{Nek12}). The representation $\rho_{f, \lambda}$ is totally odd, absolutely irreducible, and unramified outside $\mfn q$ where $q \in \mathbb{P}$ lies below $\lambda$. Moreover, for all $\mfp \in P$ satisfying $(\mfp, \mfn q) = 1$, we have
    \begin{equation}\label{eq:frob-invariants}
        \trace\bigl(\rho_{f, \lambda}(\Frob_\mfp)\bigr) = C(\mfp, f) \quad \text{and} \quad \det\bigl(\rho_{f, \lambda}(\Frob_\mfp)\bigr) = \Psi^*(\mfp)N(\mfp)^{k_0 - 1}.
    \end{equation}

    By choosing a $G_F$-stable $\mcO_{E_{f, \lambda}}$-lattice in $V_\lambda(f)$, reduction modulo $\lambda$ defines the residual Galois representation
    \begin{equation}\label{eq:residual-rep}
        \bar{\rho}_{f, \lambda} \colon G_F \longrightarrow \GL_2(\F_\lambda),
    \end{equation}
    where $\F_\lambda = \mcO_{E_f}/\lambda$. The semi-simplification of $\bar{\rho}_{f, \lambda}$ is independent of the choice of lattice. We now recall the Chebotarev density theorem (cf.~\cite{Ser81}):

    \begin{thm}[Chebotarev Density Theorem]
    \label{chebotarev_density}
        Let $G:=\bar{\rho}_{f, \lambda}(G_F)$ and let $C \subseteq G$ be a conjugacy class. Then the set
        \begin{equation*}
            \big\{\mfp \in P \mid [\bar{\rho}_{f, \lambda}(\Frob_\mfp)]_G = C\big\}
        \end{equation*}
        has Dirichlet density $\frac{|C|}{|G|}$.
    \end{thm}

\subsubsection{Inner twists of primitive forms}
    Let $f \in S_k(\mfn, \Psi)$ be a primitive Hilbert modular form over $F$.
    \begin{dfn}
        The form $f$ is said to be of CM type if there exists a nontrivial Hecke character $\Phi$ of $F$ such that $C(\mfp,f)=\Phi^{\ast}(\mfp)C(\mfp,f)$ for almost all $\mfp\in P$. We say that $f$ is non-CM if $f$ is not of CM type.
    \end{dfn}

    \begin{dfn}[Inner Twist]\label{def:inner-twist}
    Assume that the form $f$ is non-CM.
        An inner twist of $f$ is a pair $(\gamma,\Phi)$, where $\Phi$ is a Hecke character of $F$ and $\gamma\in\Aut(E_f)$, such that
        \begin{equation*}
            \gamma(C(\mfp, f)) = \Phi^{\ast}(\mfp) C(\mfp, f)
        \end{equation*}
        for almost all $\mfp\in P$.

        The automorphisms $\gamma\in\Aut(E_f)$ arising from all the inner twists $(\gamma,\Phi)$ form a subgroup $\Gamma(f)\leq \Aut(E_f).$ Moreover, $F_f=E_f^{\Gamma(f)}$. The field $F_f$ is totally real (cf.~\cite[Lemma 2.7]{KS23}), and $E_f/F_f$ is a finite Galois extension (cf.~\cite{Shi78}).
\end{dfn}

\subsection{Automorphic Representations}
    In this section, we follow the terminology of \cite{Nek12}. Let $\pi$ be an irreducible cuspidal automorphic representation of $\GL_2(\A_F)$ whose infinity type is $\pi_\infty=\bigotimes_{\nu\mid\infty}\pi_\nu$, where $\pi_\nu$ is a discrete series representation of weight $k_\nu\ge 2$, with all $k_v$ having the same parity. Let
    $\omega:\mathbb A_F^\times/F^\times\rightarrow \mathbb C^\times$
    be the central character of $\pi$. Then $\omega=|\cdot|^m\varphi$ for some $m\in\Z$ and some finite-order Hecke character $\varphi:\A_F^\times/F^\times\rightarrow \C^\times$ satisfying $\varphi_v(-1)=(-1)^{k_v}=(-1)^m$ for all $v\mid\infty$.

\subsubsection{Galois representations}
    Let $\Q(\pi)$ be the field of moduli of $\pi$ (cf.~\cite[\S B.2.2]{Nek12}). For every finite prime $\lambda$ of $\Q(\pi)$, let $V_\lambda(\pi)$ denote the two-dimensional $\Q(\pi)_\lambda$-vector space endowed with the continuous $G_F$-action attached to $\pi$ (cf.~\cite[\S B.5]{Nek12}). We write
    \begin{equation*}
        \rho_{\pi,\lambda}:G_F\longrightarrow \Aut_{\Q(\pi)_{\lambda}}\bigl(V_{\lambda}(\pi)\bigr)\simeq \GL_2\big(\Q(\pi)_{\lambda}\big)
    \end{equation*}
    for the corresponding Galois representation. We now recall the following result for the field $\Q(\pi)$.
    \begin{prop}[\cite{Nek12}, Proposition~B.2.3(4,7)]
    \label{prop:field-of-moduli}
        Let $S$ be the set of all $\nu\in P$ such that $\pi_\nu$ is not an unramified principal series. For any finite set of primes $\Sigma\supset S$, the field $\Q(\pi)$ is the compositum of the fields
        \begin{equation*}
            \Q\bigl(\lambda_\pi(\nu),\omega(\nu)\bigr)\quad \text{for } \nu\notin\Sigma,
        \end{equation*}
        where $\lambda_\pi(\nu)$ denotes the eigenvalue of the Hecke operator $T(\nu)$ acting on the spherical line of $\pi_\nu$.
    \end{prop}
    \subsubsection{Inner twists of cuspidal automorphic representations}
    Let $\pi$ be an irreducible cuspidal automorphic representation of $\GL_2(\A_F)$. We recall the following notions from Nekov{\'a}{\v r} (cf.~\cite[Proposition~B.2.7 and \S B.3]{Nek12} for more details).

\begin{dfn}
    The representation $\pi$ is said to be of CM type if there exists a totally imaginary quadratic extension $K/F$ such that $\pi \otimes \eta \simeq \pi,$ where $\eta = \eta_{K/F} \colon \A_F^\times / F^\times \to \{\pm 1\}$ is the nontrivial quadratic character associated to $K/F$. We say that $\pi$ is non-CM if $\pi$ is not of CM type.
\end{dfn}
     We now define the inner twists of $\pi$ following \cite[\S B.3]{Nek12}.

    \begin{dfn}[Inner twist of $\pi$]\label{def:inner-twist-pi}
    Assume that $\pi$ is non-CM. An \emph{inner twist} of $\pi$ is a pair $(\gamma,\chi_\gamma)$, where $\gamma\in\Aut\big(\Q(\pi)/\Q\big)$ and $\chi_\gamma:\A_F^\times/F^\times\ra\C^\times$ is a character such that
        \begin{equation*}
            {}^\gamma\pi\simeq\pi\otimes\chi_\gamma.
        \end{equation*}
        Equivalently, $\gamma\bigl(\lambda_\pi(\mfp)\bigr)=\chi_\gamma^{\ast}(\mfp)\lambda_\pi(\mfp)$ for all but finitely many $\mfp\in P$, where $\chi_\gamma^*$ is the ideal character associated with $\chi_\gamma$.

        The automorphisms $\gamma\in\Aut\big(\Q(\pi)/\Q\big)$ arising from all the inner twists $(\gamma,\chi_\gamma)$ form an abelian subgroup $\Gamma(\pi)\leq \Aut\big(\Q(\pi)/\Q\big)$.
   \end{dfn}
   We now identify the characters $\chi_\gamma$ and $\varphi$ with Galois characters $\chi_\gamma,\varphi: G_{F}\ra \Q(\pi)^{\times}$ via the reciprocity map $\mathrm{rec}_{F}:\A_F^\times/F^\times\ra G^{\mathrm{ab}}_F$ normalized by letting the uniformizers correspond to the geometric Frobenius elements. Therefore, for any subgroup $\Delta \subseteq \Gamma(\pi)$, we define
           \begin{equation*}
            F_\Delta:=\overline{F}^{\,\bigcap_{\gamma\in\Delta}\ker(\chi_\gamma)}, \, \Q(\pi)^{\Gamma(\pi)}:=\{x\in \Q(\pi)\mid\gamma(x)=x\, \forall\,\gamma\in\Gamma(\pi)\}\text{ and }F_\varphi:=\overline{F}^{\,\ker(\varphi)}.
         \end{equation*}
    We close this section by recalling the following theorem  of Nekov\'{a}\v{r}.

    \begin{thm}[\cite{Nek12}, Theorem~B.5.2]
    \label{ingredient_for_ribet_analogue}
        Let $\pi$ be non-CM and $\lambda$ be a finite prime of $\Q(\pi)$, and let $\ell\in\mathbb{P}$ be lying below $\lambda$ and set $\lambda_{\Gamma(\pi)}:=\lambda\cap \Q(\pi)^{\Gamma(\pi)}.$ Let $\Gamma(\pi)_\lambda\subseteq\Gamma(\pi)$ be the decomposition group of $\lambda$ in the extension $\Q(\pi)/\Q(\pi)^{\Gamma(\pi)}$. Then, for all but finitely many $\lambda$, there exists a basis of $V_\lambda(\pi)$ in which
        \begin{equation}\label{large_image_pi_nek}
            \rho_{\pi,\lambda} \bigl(G_{F_\varphi F_{\Gamma(\pi)_\lambda}}\bigr)=\left\{g\in\GL_2\bigl(\mcO_{\Q(\pi)^{\Gamma(\pi)},\lambda_{\Gamma(\pi)}}\bigr)\;\middle|\;
            \det(g)\in(\Z_\ell^\times)^{m-1}\right\},
        \end{equation}
        where $F_\varphi F_{\Gamma(\pi)_\lambda}$ is the compositum of $F_\varphi$ and $F_{\Gamma(\pi)_\lambda}$.
   \end{thm}

\section{Ribet's large image theorem for primitive HMFs}
    In this section, we establish an analogue of Ribet's large image theorem for the residual Galois representations attached to primitive Hilbert modular forms $f$.

\begin{thm}\label{analogue_of_Ribet}
    Let $f\in S_{k}(\mathfrak n,\Psi)$ be a non-CM primitive Hilbert modular form over $F$ of parallel even weight $k\geq 2$. Let $E_f$ and $F_f$ be as in Theorem~\ref{main_thm}, and denote their rings of integers by $\mcO_{E_f}$ and $\mcO_{F_f}$, respectively. Then for all but finitely many primes $\lambda \in \Spec(\mcO_{E_f})$, the image $\bar\rho_{f,\lambda}\bigl(G_{F_{\Psi^\ast}F_{\Gamma(f)_{\lambda}}}\bigr)$ is conjugate to
    \begin{equation*}
        \left\{g\in\GL_2(\mcO_{F_f}/\widetilde{\lambda})\;\middle|\;\det(g)\in(\F_\ell^\times)^{k-1}\right\},
    \end{equation*}
    where $\widetilde{\lambda}=\lambda\cap \mcO_{F_f}$ and $\ell\in \mathbb{P}$ lies below $\lambda$.
\end{thm}

    To establish Theorem~\ref{analogue_of_Ribet}, since $f$ has parallel even weight $k\geq2$, we consider the associated automorphic representation $\pi :=\pi(f)$ as in \cite[\S B.2.5]{Nek12}. A prerequisite for applying Nekov{\'a}{\v r}'s large image theorem, Theorem~\ref{ingredient_for_ribet_analogue}, is that the underlying representation $\pi(f)$ must be non-CM and this is the case since $f$ is non-CM.
    \begin{prop}\label{non_cm_relation}
    The automorphic representation $\pi(f)$ is non-CM.
    \end{prop}
    \begin{proof}
    Suppose in contrast that $\pi$ is of CM type. By definition, there exists a totally imaginary quadratic extension $K/F$ such that $\pi \otimes \eta \simeq \pi,$ where $\eta = \eta_{K/F} \colon \A_F^\times / F^\times \to \{\pm 1\}$ is the nontrivial quadratic character associated to $K/F$. Under this automorphic twist, for all but finitely many $\nu \in P$, we have
    $\lambda_\pi(\nu) = \eta^*(\nu) \lambda_\pi(\nu)$. Furthermore, for every prime ideal $\nu \nmid \mfn$, we have
    \begin{equation}\label{eigen_value_relation}
        C(\nu, f) = N(\nu)^{\frac{k}{2}-1}\lambda_\pi(\nu)
    \end{equation}
    by \cite[\S 12.3.8]{Nek06} and \cite[Proposition B.2.3.7]{Nek12}. Combining these relations yields
    \begin{equation*}
        C(\nu, f) = N(\nu)^{\frac{k}{2}-1}\lambda_\pi(\nu) = N(\nu)^{\frac{k}{2}-1}\eta^*(\nu)\lambda_\pi(\nu) = \eta^*(\nu)C(\nu, f)
    \end{equation*}
    for all but finitely many $\nu \in P$. Since $\eta$ is nontrivial, the ideal character $\eta^*$ is also nontrivial. Thus, $f$ must be of CM type, which is a contradiction.
\end{proof}
    Having established Proposition~\ref{non_cm_relation}, Nekov{\'a}{\v r}'s large image theorem (Theorem~\ref{ingredient_for_ribet_analogue} in the text) applies to $\pi$. Recall from \eqref{large_image_pi_nek} that for all but finitely many primes $\lambda$ of $\Q(\pi)$,
    \begin{equation}\label{large_image_pi}
        \rho_{\pi,\lambda} \bigl(G_{F_\varphi F_{\Gamma(\pi)_\lambda}}\bigr)=\left\{g\in\GL_2\bigl(\mcO_{\Q(\pi)^{\Gamma(\pi)},\lambda_{\Gamma(\pi)}}\bigr)\;\middle|\;
        \det(g)\in(\Z_\ell^\times)^{m-1}\right\}.
    \end{equation}
To prove Theorem~\ref{analogue_of_Ribet}, we must express each component of \eqref{large_image_pi} in terms of the corresponding data attached to 
$f$. We carry out this comparison in three steps.

  \begin{enumerate}
    \item[(i)] \emph{Central and nebentypus characters:} By~\cite[\S B.2.5]{Nek12}, the central character $\omega = \varphi$ of $\pi$ is of finite order, i.e., $m=0$. Let $\ell\in\mathbb P$ be the rational prime lying below $\lambda$. In Nekov{\'a}{\v r}'s normalization (cf.~\cite[\S 0.4]{Nek12} and \cite[\S 12.4.2]{Nek06}), for every $\nu\in P$ with $\nu\nmid\mfn\ell$,
                \begin{equation*}
                    \det\!\left(\mathrm{Id}-X\rho_{f,\lambda}(\Frob_\nu)\right)=1-C(\nu,f)X+\varphi(\nu)N(\nu)^{k-1}X^2,
                \end{equation*}
                whereas, by \eqref{eq:frob-invariants},
                \begin{equation*}
                    \det\!\left(\mathrm{Id}-X\rho_{f,\lambda}(\Frob_\nu)\right)=1-C(\nu,f)X+\Psi^\ast(\nu) N(\nu)^{k-1}X^2.
                \end{equation*}
                Hence, for every $\nu\in P$ with $\nu\nmid\mfn\ell$, we have
                \begin{equation}
                \label{varphi_Psi_cap_ast_relation}
                    \varphi(\nu)=\Psi^\ast(\nu).
                \end{equation}
                By Class Field Theory, these characters can be  identified with characters of $G_F$, with uniformizers corresponding to geometric Frobenius elements; thus
                $\varphi(\nu)=\varphi(\Frob_\nu)=\Psi^\ast(\Frob_\nu)=\Psi^\ast(\nu)$
                for all but finitely many $\nu\in P$. By the Chebotarev density theorem, $\varphi=\Psi^\ast$ as characters of $G_F$. Consequently,
                \begin{equation*}
                    F_\varphi=\overline F^{\,\ker(\varphi)}=\overline F^{\,\ker(\Psi^\ast)}=F_{\Psi^\ast}.
                \end{equation*}
             In addition, we have $m=0$, hence \eqref{large_image_pi} can be rewritten as
            \begin{equation}
            \label{large_image_pi_1}
                \rho_{\pi,\lambda}\bigl(G_{F_{\Psi^\ast}F_{\Gamma(\pi)_\lambda}}\bigr)=\left\{g\in\GL_2\bigl(\mcO_{\Q(\pi)^{\Gamma(\pi)},\lambda_{\Gamma(\pi)}}\bigr)\;\middle|\;\det(g)\in\Z_\ell^\times\right\}.
            \end{equation}

    \item[(ii)] \emph{Field of moduli:}
                By~\cite{Shi78}, we have $\Psi^\ast(\nu)\in E_f$, hence $\omega(\nu)\in E_f.$
                Since $k$ is even, by~\eqref{eigen_value_relation} and Proposition~\ref{prop:field-of-moduli}, we get
                \begin{equation}
                \label{same_moduli_field}
                    E_f=\Q(\pi).
                \end{equation}
                Therefore, \eqref{large_image_pi_1} can be rewritten as
                \begin{equation}\label{large_image_pi_2}
                    \rho_{\pi,\lambda}\bigl(G_{F_{\Psi^\ast}F_{\Gamma(\pi)_\lambda}}\bigr)=\left\{g\in\GL_2\bigl(\mcO_{E_f^{\Gamma(\pi)},\lambda_{\Gamma(\pi)}}\bigr)\;\middle|\;\det(g)\in\Z_\ell^\times\right\}.
                \end{equation}

    \item[(iii)] \emph{Inner twists:}
                We now show that the inner twists of $f$ and $\pi$ are equal. Since $N(\nu)^{\frac{k}{2}-1} \in \Q^\times$ is fixed under $\Aut(E_f/\Q)$, the relation in~\eqref{eigen_value_relation} gives that
                \begin{equation*}
                    \gamma\bigl(C(\nu, f)\bigr) = \Phi^\ast(\nu) C(\nu, f) \iff \gamma\bigl(\lambda_\pi(\nu)\bigr) = \Phi^\ast(\nu) \lambda_\pi(\nu)
                \end{equation*}
                for almost all $\nu \in P$, i.e., ${}^\gamma\pi \simeq \pi \otimes \Phi$. Thus, $(\gamma, \Phi)$ is an inner twist of $f$ if and only if it is an inner twist of $\pi$, i.e., $\Gamma(f) = \Gamma(\pi)$.               Hence, by~\eqref{same_moduli_field}, we have
                \begin{equation*}
                   F_f=E_{f}^{\Gamma(f)}=E_{f}^{\Gamma(\pi)}=\Q(\pi)^{\Gamma(\pi)}.
               \end{equation*}
               Therefore, \eqref{large_image_pi_2} can be rewritten as
                \begin{equation}\label{large_image_pi_3}
                    \rho_{\pi,\lambda}\bigl(G_{F_{\Psi^\ast}F_{\Gamma(f)_\lambda}}\bigr)=\left\{g\in\GL_2\bigl(\mcO_{F_f,\lambda_{\Gamma(f)}}\bigr)\;\middle|\;\det(g)\in\Z_\ell^\times\right\}.
                \end{equation}
    \end{enumerate}
    In~\eqref{large_image_pi_3}, the RHS is expressed entirely in terms of the data attached to $f$, except that the LHS involves the automorphic Galois representation $\rho_{\pi,\lambda}$. It remains to compare the representation $\rho_{\pi,\lambda}$ with the Galois representation $\rho_{f,\lambda}$ attached to $f$. By~\cite[\S B.2.5]{Nek12}, for every $\lambda \in\Spec(\mcO_{E_f})$, we have
    \begin{equation*}
        V_{\lambda}(\pi)=V_{\lambda}(f)\left(r\right),
    \end{equation*}
    where $r:=\frac{k}{2}-1$. By the definition of the Tate twist, the corresponding $\lambda$-adic Galois representations satisfy
    \begin{equation}
    \label{image_relation_GR}
        \rho_{\pi,\lambda}(\sigma)=\chi_{\mathrm{cyc},\ell}(\sigma)^r\rho_{f,\lambda}(\sigma) \text{ for } \sigma\in G_F,
    \end{equation}
    where $\ell \in \mathbb{P}$ lying below $\lambda$ and $\chi_{\mathrm{cyc},\ell}:G_F\to\Z_\ell^\times$ denotes the $\ell$-adic cyclotomic character.

We now combine \eqref{large_image_pi_3} and \eqref{image_relation_GR}. The proof of Theorem~\ref{analogue_of_Ribet} relies on the following two lemmas. For notational simplicity, set $G_{\lambda}:=G_{F_{\Psi^\ast}F_{\Gamma(f)_{\lambda}}}.$

    \begin{lem}
    \label{determinant_image_rho_f}
        Let $f$ be as in Theorem~\ref{analogue_of_Ribet}. For almost all $\lambda$, $\det\bar\rho_{f,\lambda}(G_{\lambda})=(\F_\ell^\times)^{k-1}.$
    \end{lem}

    \begin{proof}
         By \eqref{large_image_pi_3}, for almost all $\lambda$,
        $\rho_{\pi,\lambda}(G_{\lambda})=\left\{g\in\GL_2(\mcO_{F_f,\widetilde{\lambda}})
            \;\middle|\;\det(g)\in\Z_\ell^\times\right\},$
    where $\widetilde{\lambda}=\lambda\cap \mcO_{F_f}$. Now, reducing modulo $\lambda$, the equality becomes
    \begin{equation}
    \label{image_of_bar_rho_pi}
        \bar{\rho}_{\pi,\lambda}(G_{\lambda})=\left\{g\in\GL_2(\mcO_{F_f}/\widetilde{\lambda})\;\middle|\;\det(g)\in\F_\ell^\times \right\}.
    \end{equation}

    \begin{itemize}
         \item By~\cite[\S B.2.5]{Nek12} and \eqref{varphi_Psi_cap_ast_relation}, we have $\det\bar{\rho}_{\pi,\lambda}=\overline{\chi}_{\mathrm{cyc},\ell}^{-1}\overline{\Psi}^\ast.$ Since $G_{\lambda}\subseteq G_{F_{\Psi^\ast}}=\ker(\Psi^\ast),$ the character $\Psi^\ast$ is trivial on $G_{\lambda}$. So, we obtain
    \begin{equation}
    \label{det_rho_pi_chi}
        \det\bar\rho_{\pi,\lambda}(\sigma)=\overline{\chi}_{\mathrm{cyc},\ell}(\sigma)^{-1} \text{ for }\sigma\in G_{\lambda}.
    \end{equation}
    Thus,  we have $\overline{\chi}_{\mathrm{cyc},\ell}(G_{\lambda})=\F_\ell^\times$, since $\det\bar\rho_{\pi,\lambda}(G_{\lambda})=\F_\ell^\times$  (cf. \eqref{image_of_bar_rho_pi}).

    \item   Now,  reducing modulo $\lambda$,~\eqref{image_relation_GR} becomes $\bar\rho_{\pi,\lambda}(\sigma)=\overline{\chi}_{\mathrm{cyc},\ell}(\sigma)^r\bar\rho_{f,\lambda}(\sigma)$ for $\sigma\in G_F$. In particular, we have
    \begin{equation*}
        \det\bar\rho_{f,\lambda}(\sigma)=\overline{\chi}_{\mathrm{cyc},\ell}(\sigma)^{-2r}\det\bar\rho_{\pi,\lambda}(\sigma) \text{ for }  \sigma\in G_F.
    \end{equation*}
    \end{itemize}
    By~\eqref{det_rho_pi_chi} for $\sigma\in G_{\lambda}$, the last equality reduces to $\det\bar\rho_{f,\lambda}(\sigma)=\overline{\chi}_{\mathrm{cyc},\ell}(\sigma)^{-(2r+1)}= \overline{\chi}_{\mathrm{cyc},\ell}(\sigma)^{-(k-1)}$. Since $\overline{\chi}_{\mathrm{cyc},\ell}(G_{\lambda})=\F_\ell^\times$, we have $\det\bar\rho_{f,\lambda}(G_{\lambda})= (\F_\ell^\times)^{k-1}$.
\end{proof}

    \begin{lem}
    \label{SL2_contained_rho_f}
         Let $f$ be as in Theorem~\ref{analogue_of_Ribet}. For almost all $\lambda$, $\SL_2(\mcO_{F_f}/\widetilde{\lambda})\subseteq\bar\rho_{f,\lambda}(G_{\lambda}).$
    \end{lem}

   \begin{proof}
    By~\eqref{image_of_bar_rho_pi}, the inclusion  $\SL_2(\mcO_{F_f}/\widetilde{\lambda}) \subseteq \bar\rho_{\pi,\lambda}(G_{\lambda})$ gives us the inclusion of commutator subgroups $[\SL_2(\mcO_{F_f}/\widetilde{\lambda}), \SL_2(\mcO_{F_f}/\widetilde{\lambda})] \subseteq [\bar\rho_{\pi,\lambda}(G_{\lambda}), \bar\rho_{\pi,\lambda}(G_{\lambda})]$. If $|\mcO_{F_f}/\widetilde{\lambda}|> 3$, the group $\SL_2(\mcO_{F_f}/\widetilde{\lambda})$ is perfect, i.e., the commutator subgroup of $\SL_2(\mcO_{F_f}/\widetilde{\lambda})$ is itself. So, after excluding finitely many $\lambda$, we get $\SL_2(\mcO_{F_f}/\widetilde{\lambda})\subseteq [\bar\rho_{\pi,\lambda}(G_{\lambda}), \bar\rho_{\pi,\lambda}(G_{\lambda})]$. Since the residual representations $\bar\rho_{\pi,\lambda}$ and $\bar\rho_{f,\lambda}$ differ only by a scalar, we get
\begin{equation*}
    \SL_2(\mcO_{F_f}/\widetilde{\lambda}) \subseteq [\bar\rho_{\pi,\lambda}(G_{\lambda}), \bar\rho_{\pi,\lambda}(G_{\lambda})] = [\bar\rho_{f,\lambda}(G_{\lambda}), \bar\rho_{f,\lambda}(G_{\lambda})] \subseteq \bar\rho_{f,\lambda}(G_{\lambda}).
\end{equation*}
\end{proof}

  We are now ready to prove Theorem~\ref{analogue_of_Ribet}.

    \begin{proof}[Proof of Theorem~\ref{analogue_of_Ribet}]
    Let
        $H =\bar\rho_{f,\lambda}(G_{\lambda})$
    and $$ G =\left\{g\in\GL_2(\mcO_{F_f}/\widetilde{\lambda}) \;\middle|\;\det(g)\in(\F_\ell^\times)^{k-1}\right\}.$$
    We aim to prove that $H=G$. Since $\bar\rho_{\pi,\lambda}(\sigma)= \overline{\chi}_{\mathrm{cyc},\ell}(\sigma)^r \bar\rho_{f,\lambda}(\sigma)$ for $\sigma\in G_F$, we get that  $H \subseteq \GL_2(\mcO_{F_f}/\widetilde{\lambda})$, by~\eqref{image_of_bar_rho_pi}. By Lemma~\ref{determinant_image_rho_f}, we have $\det(H) = (\F_\ell^\times)^{k-1}$, hence $H \subseteq G$.

 For the reverse inclusion, let $g\in G$. By Lemma~\ref{determinant_image_rho_f}, $\det(g)\in (\F_\ell^\times)^{k-1}= \det(H)$, there exists $h\in H$ such that $\det(g)=\det(h)$. By Lemma~\ref{SL2_contained_rho_f}, we get
    $h^{-1}g\in\SL_2(\mcO_{F_f}/\widetilde{\lambda})\subseteq H$ and hence $g\in H$. This proves the claim.
\end{proof}

    We now recall a group theoretic estimate of Koo, Stein, and Wiese.
\begin{prop}[\cite{KSW08}, Proposition~2.1(c)]
\label{density_characteristic_poly}
    Let $q$ be a prime power, $r\in\N$, and let $R\subseteq \widetilde{R}\subseteq \F_{q^r}^{\times}$ be subgroups. Let $G\subseteq\left\{g\in\GL_2(\F_{q^r}) :\det(g)\in\widetilde{R}\right\}$ be a subgroup such that $H:=\left\{g\in\GL_2(\F_q):\det(g)\in R\right\}$ is a normal subgroup of $G$. Let $P(X)=X^2-aX+b\in\F_{q^r}[X],$ then $\sum_C |C|\leq2\lvert\widetilde{R}/R\rvert\,(q^2+q),$ where the sum runs over all conjugacy classes $C$ of $G$ whose characteristic polynomial is $P(X)$.
\end{prop}

    We now combine the preceding large image theorem, i.e.,  Theorem~\ref{analogue_of_Ribet}, with the above group theoretic estimate to obtain the following:
\begin{thm}
    Let $f, F_f$ be as in Theorem~\ref{main_thm} and $L\subset F_{f}$ be any proper subfield. Then the set
    \begin{equation*}
        \left\{ \mfp \in P \;\middle|\; \frac{C(\mfp,f)^2}{\Psi^{\ast}(\mfp)}\in L \right\}
    \end{equation*}
    has Dirichlet density $0$.
\end{thm}
\begin{proof}

    Let $\mcO_{L}$ be the ring of integers of $L$. We define
    \begin{equation*}
        \mcS := \left\{ \widetilde{\lambda} \in \Spec(\mcO_{F_f}) \;\middle|\; [\mcO_{F_f}/\widetilde{\lambda} : \mcO_L/(\widetilde{\lambda} \cap \mcO_L)] \geq 2 \right\}.
    \end{equation*}
    By the Chebotarev density theorem, $\mcS$ is infinite. After excluding finitely many primes from $\mcS$, we may assume without loss of generality that Theorem~\ref{analogue_of_Ribet} holds for every prime $\lambda \in \Spec(\mcO_{E_f})$ lying above any $\widetilde{\lambda} \in \mcS$.

    Fix $\widetilde{\lambda} \in \mcS$. Let $\ell\in \mathbb{P}$ lie below $\widetilde{\lambda}$, and let $\lambda \in \Spec(\mcO_{E_f})$ be any prime lying over $\widetilde{\lambda}$. Set
    \begin{equation*}
        \F_q := \mcO_L/(\widetilde{\lambda} \cap \mcO_L), \quad \F_{q^r} := \mcO_{F_f}/\widetilde{\lambda}, \quad \text{and} \quad \F_{q^{rs}} := \mcO_{E_f}/\lambda.
    \end{equation*}
    By definition of $\mcS$, we have $r \geq 2$. Let $R:=(\F_{\ell}^{\times})^{(k-1)}$ and  $W\subseteq \F_{q^{rs}}^{\times}$ denote the image of $\Psi^{\ast}$ modulo $\lambda$.  Let $\widetilde{R}:=\langle R,W\rangle$ be the subgroup of $\F_{q^{rs}}^{\times}$ generated by $R$ and $W$. Note that $|W|\leq \left|(\mcO_F/\mfn)^{\times}\right|$ and $|\widetilde{R}|\leq |R|\,|W|$.

    Let $H:=\{g\in \GL_{2}(\F_{q^r}) \mid \det(g)\in (\F_{\ell}^{\times})^{(k-1)}\}$ and $G:=\Gal(\overline{F}^{\ker(\bar{\rho}_{f,\lambda})}/F)$, which we can identify with the image of $\bar{\rho}_{f,\lambda}$. By Theorem~\ref{analogue_of_Ribet} and ~\eqref{eq:frob-invariants}, \eqref{eq:residual-rep}, we get
    \begin{equation*}
        H\subseteq G\subseteq \{g\in \GL_{2}(\F_{q^{rs}}) \mid \det(g)\in \widetilde{R}\},
    \end{equation*}
    with $H$ being normal in $G$, since $G_{\lambda}$ is normal in $G_F$. Now consider the set
    \begin{equation*}
        M_{\widetilde{\lambda}}:=\bigsqcup_{C}\left\{\mfp\in P\;\middle|\; [\bar{\rho}_{f,\lambda}(\Frob_{\mfp})]_{G}=C\right\}\supseteq \left\{\mfp\in P\;\middle|\;\overline{\frac{C(\mfp,f)^2}{\Psi^{\ast}(\mfp)}}\in \F_{q} \right\},
    \end{equation*}
    where the reduction modulo $\widetilde{\lambda}$ of an element $x\in \mcO_{F_{f}}$ is denoted by $\bar{x}$ and $C$ runs over the conjugacy classes of $G$ with characteristic polynomials equal to some $x^2 - ax + b \in \F_{q^{rs}}[x]$ such that
    \begin{equation*}
        a^2\in\left\{t\in\F_{q^{rs}}\;\middle|\;\exists u\in \F_q, \exists w\in W \text{ such that }t=uw\right\}.
    \end{equation*}
   Note that, $b\in \widetilde{R}$. There are at most $2q|W|^2\,|R|$ such polynomials. By Theorem~\ref{chebotarev_density}, we have $\delta_{D}(M_{\widetilde{\lambda}})=\sum_{C}\frac{|C|}{|G|}$.  Now, by Proposition~\ref{density_characteristic_poly} and $r \geq 2$, we get
    \begin{align*}
        \delta_{D}(M_{\widetilde{\lambda}})&\leq \frac{2q|W|^2\,|R|\times 2 |\widetilde{R}/R|(q^{2r} + q^r)}{|G|}\\
        &\leq\frac{4|W|^3\,|R|q(q^{2r}+q^r)}{|R|\,|\SL_{2}(\F_{q^r})|}\\
        &\leq\frac{4|W|^3q(q^{2r}+q^r)}{q^{3r}-q^r}=O\left(\frac{1}{q^{r-1}}\right)\leq O\left(\frac{1}{q}\right),
    \end{align*}
    where we used the bound $|G| \geq |H| = |R|\,|\SL_2(\F_{q^r})|$ with $|\SL_2(\F_{q^r})| = q^{3r} - q^r$, and the estimate $|\widetilde{R}/R| \leq |W|$.  Now, set
    \begin{equation*}
        A:=\left\{\mfp\in P\;\middle|\;\frac{C(\mfp,f)^2}{\Psi^\ast(\mfp)}\in L\right\}.
    \end{equation*}
    Since $C(\mfp,f)$ is an algebraic integer and $\Psi^\ast(\mfp)$ is a root of unity, we get $\frac{C(\mfp,f)^2}{\Psi^\ast(\mfp)} \in L\cap \mcO_{E_f}=\mcO_L.$
    Therefore, for every $\widetilde{\lambda}\in\mcS$, we have $A\subseteq M_{\widetilde{\lambda}}$ and hence, the upper Dirichlet density of $A$ is bounded above by  $\delta_D(M_{\widetilde{\lambda}})\leq O\!\left(\frac1q\right).$ Since $q$ is unbounded as $\widetilde{\lambda}$ varies in $\mcS$, it follows that the upper Dirichlet density of $A$ is $0$. Consequently, the set $A$ has Dirichlet density $0$.
\end{proof}
    We are now in a position to give a proof of the Main Theorem (cf. Theorem~\ref{main_thm} in the text)  of this article.
\begin{proof}[Proof of Theorem~\ref{main_thm}]
    Since $E_f/\Q$ is a finite and separable extension (cf.~\cite{Shi78}), the field $F_f$ has only finitely many subfields. If $\Q\big(\frac{C(\mfp,f)^2}{\Psi^\ast(\mfp)}\big)\subsetneq F_f,$ then $\frac{C(\mfp,f)^2}{\Psi^\ast(\mfp)}$ lies in some proper subfield $L\subset F_f$. By the preceding theorem, for each such $L$, the set of prime ideals $\mfp$ satisfying this condition has Dirichlet density $0$. Since $F_f$ has only finitely many proper subfields, the union of these density zero sets also has Dirichlet density $0$. Therefore,
    \begin{equation}
    \label{density_F_f}
         \delta_D\left(\left\{\mfp\in P\;\middle|\;\Q\left(\frac{C(\mfp,f)^2}{\Psi^\ast(\mfp)}\right)=F_f\right\}\right)=1\end{equation}
    We now prove the second part of the theorem, regarding the single generation of $E_f$.
    By~\cite[Lemma~2.7]{KS23}, we get $\frac{C(\mfp,f)^2}{\Psi^\ast(\mfp)}\in \Q(C(\mfp,f))$ for all $\mfp\in P$ with $\mfp\nmid \mfn$. If $f$ has no nontrivial inner twists, then $E_f=F_f$. Then, by~\eqref{density_F_f} we get
    \begin{equation*}
        \delta_D\left(\left\{\mfp\in P\;\middle|\;\Q(C(\mfp,f))=E_f\right\}\right)=1.
    \end{equation*}
\end{proof}

    The following corollaries are immediate consequences of Theorem~\ref{main_thm}, by arguments analogous to those in \cite{KS23}, we therefore omit the proofs.

    For any subgroup $\Gamma'\leq\Gamma :=\Gamma(f)$, let $\{\Psi_{\gamma}^{\ast}\mid\gamma\in\Gamma'\}$ denote the corresponding inner twist characters, regarded as characters of $G_F$, and set
    \begin{equation*}
         H_{\Gamma'}:=\bigcap_{\gamma\in\Gamma'}\ker(\Psi_\gamma^\ast)\text{ and } K_{\Gamma'}:=\overline F^{\,H_{\Gamma'}}.
    \end{equation*}
    \begin{cor}\label{density_membership_intermediate}
        Let $f$, $E_f$, and $F_f$ be as in Theorem~\ref{main_thm}, and let $K\subseteq E_f$ be a subfield. Then the following hold:
        \begin{enumerate}
            \item[(i)] If $F_f\nsubseteq K$, then
            \begin{equation*}
                \delta_D\left(\left\{\mfp \in P \;\middle|\; C(\mfp, f) \in K\right\}\right) = \delta_D\left(\left\{\mfp \in P \;\middle|\; \Q(C(\mfp, f)) = K\right\}\right) = 0.
            \end{equation*}
            \item[(ii)] If $F_f \subseteq K$, let $\Gamma' \leq \Gamma$ be the unique subgroup such that $K = E_f^{\Gamma'}$. Then
            \begin{equation*}
                \delta_D\left(\left\{\mfp \in P \;\middle|\; C(\mfp, f) \in K\right\}\right) = \frac{1}{[K_{\Gamma'} : F]}.
            \end{equation*}
        \end{enumerate}
    \end{cor}

    \begin{cor}\label{density_generation_intermediate}
        Let $f$ and $F_f$ be as in Theorem~\ref{main_thm}, and let $K$ be an intermediate subfield such that $F_f \subseteq K \subseteq E_f$. Let $\Gamma' \leq \Gamma$ be the unique subgroup such that $K = E_f^{\Gamma'}$. Then the set
        \begin{equation*}
            \left\{\mfp \in P \;\middle|\; \Q(C(\mfp, f)) = K\right\}
        \end{equation*}
        has Dirichlet density equal to that of the set
        \begin{equation*}
            \left\{\mfp \in P \;\middle|\; \Psi_\gamma^\ast(\mfp) = 1 \text{ for all } \gamma \in \Gamma'\,\text{ and }\,\Psi_\omega^\ast(\mfp) \neq 1 \text{ for all } \omega \in \Gamma \setminus \Gamma'\right\}.
        \end{equation*}

    \end{cor}

\section*{Acknowledgements}

The authors are grateful to Prof. Ariel Weiss for bringing our attention to the work of Nekov\'{a}\v{r}~\cite[\S B]{Nek12}.

\bibliographystyle{plain, abbrv}
	
\end{document}